\documentclass[12pt]{article}
\usepackage[cp1251]{inputenc}
\usepackage{epic}
\usepackage{eepic}
\usepackage{graphics}
\usepackage{amsfonts}
  \usepackage{amsthm}
 \usepackage{amsmath}
 \usepackage{amscd}
 \usepackage{amssymb}
 \usepackage{latexsym}
\usepackage{graphicx}

 \numberwithin{equation}{section}
 \newtheorem{thm}{Theorem}
 \newtheorem{lm}{Lemma}
 \newtheorem{prop}{Proposition}
 \theoremstyle{definition}
 
\newcommand{\ds}{\displaystyle}

\title{  Non-Integrability of the Trapped Ionic System II}
\author{ Georgi Georgiev \\
 Faculty of Mathematics and Informatics, \\
 Sofia University ``St. Kl. Ohridski'', \\
 5 James Bourchier Blvd., 1164 Sofia, Bulgaria}
\date{}

\begin{document}

\maketitle

\begin{abstract}
The idea for writing this paper is to analyze the differences in the results in \cite{G} and \cite{Comment I}, correct the errors in \cite{G}, and remove the inappropriate comments in \cite{G} about the results in \cite{Trapped}. It is necessary and should be noted that the inaccuracies \cite{Trapped} have been removed.

\end{abstract}
{\bf Key words:} Hamiltonian system, Meromorphic non-integrability, Variational equation, Monodromy group, Differential Galois group
\bigskip

\section{Introduction}

This paper is a continuation of \cite{G} and for this reason we will omit the introduction and research data, referring to the original article (see also \cite{Trapped} and \cite{TrappedI} for details). I will not go over the classical and modern theory of meromorphic non-integrability studies, I will simply refer to \cite{Lyap}, \cite{MSim}, \cite{MR1},\cite{MR2}, \cite{Jarno1},\cite{Jarno}, \cite{DarbouxPoint}, \cite{DarbouxPoint1}, \cite{SzWo}, \cite{Chur} and \cite{Berger}.
 We study two dimensional model 
\begin{equation}
\label{1.1}
H=\frac{1}{2}(p_r^2+p_z^2)+Ar^2+Bz^2+Cz^3+Dr^2z+Ez^4+Fr^2z^2+Gr^4,
\end{equation}
 where $A$, $B$, $C$, $D$, $E$, $F$, and $G$  are  a appropriate
real constants for existing an additional  meromorphic integral of motion.

Let  we denote with
$p:=\sqrt{1+\frac{4F}{E}}$, and with $z_i$, $i=1,\, \dots , 4$ the roots of polynomial $Ez^4+Cz^3+Bz^2+h=0$, where $h$ is a constant (we suppose that $z_i\ne z_j$, for $i\ne j$). 

The result of this article is as follows:
\begin{thm}
\label{main2D_1}

a)\, Assume that  $p\notin\mathbb{Q}$, then the system (\ref{1.1})
has no an additional analytic first integral;

b)\, Let $p\, \in   \mathbb{Q}$ if  $N(p)\ge 4$, then the system (\ref{1.1})
has no an additional meromorphic first integral;

c)\, Let $p\, \in   \mathbb{Q}$ and   $N(p)\le 3$, 
then the system (\ref{1.1}) has no an additional meromorphic first integral, if at least one of the following  conditions is true:

c.0)\, if $V_0=0$, ($E=F=G=0$), then  \\
c.01)\, $D\ne 0$,  $\frac{C}{D}\notin \{\frac{1}{3},\,2,\,\frac{16}{3}\}$, or;\\  
c.02)\, $D\ne 0$, $\frac{C}{D}=\frac{1}{3}$, $A\ne B$;\\
c.03)\, $D\ne 0$, $\frac{C}{D}=\frac{16}{3}$, and $16A\ne 5B$;

c.1)\, if $V_1=z^4$, ($F=G=0$), then $D\ne 0$;

c.2)\, if $V_1=r^4$, ($E=F=0$), then  $D\ne 0$;

c.3)\, if $V_{2}= (r^2+z^2)^2$, ($E=G=\ds{\frac{F}{2}}$), then \\
c.31)\, if $D\ne0$, $\frac{C}{D}\notin \{\frac{1}{3},\,2,\,\frac{16}{3}\}$, or;\\
c.32)\, if $D\ne0$, $\frac{C}{D}=\frac{1}{3}$, then $A\ne B$;\\
c.33)\,  if $D\ne0$, $\frac{C}{D}=\frac{1}{3}$, $A= B$, then $A\ne 0$ and $4A\ne 15D^2$;\\
c.34)\, if $D\ne0$, $\frac{C}{D}=2$, then $D^2\ne B- 2A$;\\
c.35)\, if $D\ne0$, $\frac{C}{D}=\frac{16}{3}$, then $16A\ne5B$; \\
c.36)\, if $D\ne0$, $\frac{C}{D}=\frac{16}{3}$, $16A=5B$, then $20A\ne -63D^2$;\\

c.4)\, if $V_3=r^4+6r^2z^2+z^4$, ($G=E=6F$), then \\
c.41)\, if $D\ne0$, $\frac{C}{D}\notin \{\frac{1}{3},\,2,\,\frac{16}{3}\}$, or;\\
c.42)\, if $D\ne0$, $\frac{C}{D}=\frac{1}{3}$, then $A\ne B$;\\
c.43)\,  if $D\ne0$, $\frac{C}{D}=2$, then $D^2\ne B- 2A$;\\
c.44)\, if $D\ne0$, $\frac{C}{D}=\frac{16}{3}$, then $16A\ne 5B$; \\
c.45)\, if $D\ne0$, $\frac{C}{D}=\frac{16}{3}$, $16A=5B$, then $20A\ne -63D^2$;\\

c.5)\,  if $V_4=z^4+\frac{1}{4}\alpha r^4$, ($\frac{E}{G}=\frac{4}{\alpha}$, $\alpha\ne 0$), then $D\ne 0$;

c.6)\, if $V_4=r^4+\frac{1}{4}\alpha z^4$, ($\frac{E}{G}=\frac{\alpha}{4}$, $\alpha\ne 0$), then $D\ne 0$;

c.7)\, if $V_5=16r^4+3r^2z^2+z^4$, ($G=16E$, $3G=4F$), then ;\\
c.71)\, if $D\ne0$, $\frac{C}{D}\notin \{\frac{1}{3},\,\frac{16}{3}\}$, or;\\
c.72)\, if $D\ne0$, $\frac{C}{D}=\frac{1}{3}$, then $A\ne B$;\\
c.73)\, if $D\ne0$, $\frac{C}{D}=\frac{16}{3}$, $16A\ne 5B$;\\
c.74)\, if $D\ne0$, $\frac{C}{D}=\frac{16}{3}$, $16A=5B$, then $16A\ne -15D^2$;\\
c.75)\, if $D=0$ and $C=0$,  then  $A\ne 4B$;\\

c.8)\, if $V_5=16z^4+3r^2z^2+r^4$, ($E=16G$, $3E=4F$), then\\
c.81)\, if $D\ne0$, $\frac{C}{D}\notin \{\frac{1}{3},\,2,\,\frac{16}{3}\}$, or;\\
c.82)\, if $D\ne0$, $\frac{C}{D}=\frac{1}{3}$, then $A\ne B$;\\
c.83)\,  if $D\ne0$, $\frac{C}{D}=\frac{1}{3}$, $A= B$, then  $8A\ne 5D^2$;\\
c.84)\,  if $D\ne0$, $\frac{C}{D}=2$, then $ B\ne 4A$;\\
c.85)\, if $D\ne0$, $\frac{C}{D}=\frac{16}{3}$, $16A\ne 5B$;\\
c.86)\, if $D\ne0$, $\frac{C}{D}=\frac{16}{3}$, $16A=5B$, then $16A\ne -15D^2$;\\
c.87)\, if $D=0$ and $C=0$,  then $A\ne 0$ and $A\ne 3B$;\\

c.9)\, if $V_6=z^4+6r^2z^2+8r^4$ ($E=8G$, $E=6F$), then\\
c.91)\, if $D\ne0$,  or;\\
c.92)\,  if $D=0$ and $C=0$,   $A\ne 4B$;\\

c.10)\, if $V_6=r^4+6r^2z^2+8z^4$ ($F=6G$, $3E=4F$), then\\
c.101)\, if $D\ne0$,  or;\\
c.102)\,  if $D=0$ and $C=0$,   $4A\ne B$.
 \end{thm}
Here $ N (r) $ is the positive denominator of the irreducible $ r \in \mathbb {Q} $, with $V_i$, $i=1,\dots,6$ denotes all integrable homogeneous potentials of the fourth degree ($V_0$ is when there is no such potential). It should be noted that the notations in this paper are not quite the same as in \cite{DarbouxPoint}. The study of exactly these specific integrable potentials does not restrict the community of the problem applying {\bf Remark 1} in paper \cite{DarbouxPoint1}. 
This is probably the longest theorem that I've ever formulated, unfortunately this is the result of researching this problem. The idea of the proof is to study for integrability the sums of the potentials of the second degree (in our case it is integrable) with all the integrable potentials of the third degree, then to the Integrable sums are added to all possible integrable potentials of the fourth degree and again check for integrability. The homogeneous integrable potentials of degree three and four have been studied exhaustively  in \cite{Jarno1}, \cite{Jarno}, \cite{NakYo}, \cite{DarbouxPoint1}, \cite{DarbouxPoint}, \cite{ME}, \cite{Yo}.

\section{$VE_1$}

 The Hamiltonian equations are:
\begin{align}
\label{1.2}
\dot r  &=  p_r,\, \dot p_r  =  -(2Ar+2Drz+2Frz^2+4Gr^3),\nonumber  \\
\dot z  &=  p_z,\, \dot p_z  =  -(2Bz+3Cz^2+Dr^2+4Ez^3+2Fr^2z).
\end{align}
for existing an additional  integral of motion (here as usual $\dot{}=\frac{d}{dt}$).

In the our study we will suppose that $(F,\, D)\ne (0,\, 0)$, because if we assume the opposite, the variables in the considered system are separated, i.e. the system is integrable.

First, we find a partial solution for (\ref{1.2}). Let we put
$r = p_r = 0$ in (\ref{1.2}) and we have
$$
\ddot z=-(2Bz+3Cz^2+4Ez^3),
$$
multiplying by $\dot z$ and integrating by the time $t$ we have
\begin{equation}
\label{1.3}
{\dot z}^2=-2(Ez^4+Cz^3+Bz^2+h),
\end{equation}
where $h$ is a  constant. Further, we follow the procedures for
Ziglin-Morales-Ramis theory and we find an invariant manifold
$(r,\,p_r,\,z,\, p_z)=(0,\,0,\, z,\, \dot z)$ here $z$ is the
solution of (\ref{1.3}). According to theory, the solution of
(\ref{1.3}) must be a rational function of Weierstrass
$\wp$-function, but it is not important for us right now. Finding
the Variation Equations (VE) we have $\xi_{11}=dr$, $\eta_{11}=dp_r$,
$\xi_{12}=dz$, and  $\eta_{12}=dp_z$ and we obtain:
\begin{align}
\label{VE1.4}
 \ddot \xi_{11 }&=  -2(A+Dz+Fz^2)\xi_{11},\nonumber\\
\ddot \xi_{12}  &=  -2(B+3Cz+6Ez^2)\xi_{12}.
\end{align}
Next we change  the variable in equations (\ref{VE1.4}) by
$\xi_{1i}(t)=\xi_{1i}(z(t))$, where $z(t)$ is a solution of (\ref{1.3}) and
we have $\frac{d\xi_{1i}(t)}{dt}=\frac{d\xi_{1i}}{dz}. \frac{dz(t)}{dt}$,
for $ i=1,\, 2$  and
\begin{eqnarray*}
\frac{d^2\xi_{11}}{dt^2} &=&\frac{d^2\xi_{11}}{dz^2}\left (\frac{dz}{dt}\right )^2+\frac{d\xi_{11}}{dz}\frac{d^2z}{dt^2}\\
&=& -2(Ez^4+Cz^3+Bz^2+h)\frac{d^2\xi_{11}}{dz^2}-(2Bz+3Cz^2+4Ez^3)\frac{d\xi_{11}}{dz}\\
&+& 2(A+Dz+Fz^2)\xi_{11}=0,
\end{eqnarray*}
and
\begin{eqnarray*}
\frac{d^2\xi_{12}}{dt^2} &=&\frac{d^2\xi_{12}}{dz^2}\left (\frac{dz}{dt}\right )^2+\frac{d\xi_{12}}{dz}\frac{d^2z}{dt^2}\\
&=& -2(Ez^4+Cz^3+Bz^2+h)\frac{d^2\xi_{12}}{dz^2}-(2Bz+3Cz^2+4Ez^3)\frac{d\xi_{12}}{dz}\\
&+& 2(B+3Cz+6Ez^2)\xi_{12}=0.
\end{eqnarray*}
If we denote with  $' =\frac{d}{dz}$ we obtain for the VE two
Fuchsian linear differential equations with four singularities:
\begin{eqnarray}
\label{VE1.5}
 \xi_{11}''  +\frac{4Ez^3+3Cz^2+2Bz}{2(Ez^4+Cz^3+Bz^2+h)}\xi_{11}'  -\frac{A+Dz+Fz^2}{Ez^4+Cz^3+Bz^2+h}\xi_{11}=0, \nonumber \\
 \xi_{12}''  +\frac{4Ez^3+3Cz^2+2Bz}{2(Ez^4+Cz^3+Bz^2+h)}\xi_{12}'
-\frac{B+3Cz+6Ez^2}{Ez^4+Cz^3+Bz^2+h}\xi_{12}=0.
\end{eqnarray}

 The equations ( \ref{VE1.5}) are a Fuchsian and have five regular singularities $z_i$, for $i=1,\, \dots , 4$ and $z=\infty$.

Here we study the equations ( \ref{VE1.5}) for a Liouvillian
solutions. Existence of such solutions of a linear differential equation is equivalent
to the solvability of the identity component of its Galois group.

First we find conditions for branching for
solutions of this equation. We will use Lyapunov's  idea to prove
non-integrability using the branching of the solutions of the
variational equations around an appropriate partial solution. If the
solutions of the VE branching then the Hamiltonian system has not
additional first integral. (See \cite{Lyap} for details.) In this section we will assume that $B\ne 0$ and $E\ne 0$. We find
the  indicial equations to the singular points $z_i$ and $z=\infty$.
We have $\lambda^2-\frac{1}{2}\lambda=0$, roots area
$\lambda_1=0$, $\lambda_2=\frac{1}{2}$,  for $z_i$, second one is
$\rho^2-\rho -\frac{F}{E}=0$, with roots
$\rho_k=\frac{1\pm\sqrt{1+\frac{4F}{E}}}{2}$, $k=1,2$ for
$z=\infty$. The solutions around $z_i$ have no branching  and near $z=\infty$
we have branching - when $\rho_k, \, \notin \mathbb{Q}, \, k=1,\,2
$. This is an application of the Frobenius's method: the fundamental system of solutions around of the singular
points are presented in the form
 $\Phi_j(z)=z^{\lambda_j}\Psi_j(z)$, where $\Psi_j(z)$ are holomorphic functions (locally) for $j=1,2$. Then the solutions of ( \ref{VE1.5}
 are branching for $\lambda_j, \, \notin \mathbb{Q}, \, j=1,\, 2 $. For $z=\infty$   we do the same and obtain that $\rho_k, \, \notin \mathbb{Q}, \, k=1,\,2 $. It should be noted that it is possible, if $\lambda_j\in \mathbb{Q}$ and $\rho_k\in \mathbb{Q}$, we also can achieve not branching by a standard changing of the variables. This change of variables does not affect the commutativity of the unity component of the Galois group.
We have proved the following proposition.
\begin{prop}
\label{Lyap}
Let  $p=\sqrt{1+4\frac{F}{E}}\notin\mathbb{Q}$, then the system (\ref{1.2})
has no additional holomorphic first integral.
\end{prop}
This proves $a)$ in Theorem \ref{main2D_1}.

Now we use the notation
$F=\ds{\frac{p^2-1}{4}}E$, where
$p\, \in \mathbb{Q}$. For the first of eq. ( \ref{VE1.5}) we put
$$c_1(z):=\frac{4Ez^3+3Cz^2+2Bz}{2(Ez^4+Cz^3+Bz^2+h)}, $$
and
$$c_2(z):=-\frac{A+Dz+\frac{p^2-1}{4}Ez^2}{Ez^4+Cz^3+Bz^2+h}.$$
 
In the terminology of \cite{Chur} we get 
$a_{\infty}=\lim_{z\to\infty} z c_1(z)=2$, $b_{\infty}=\lim_{z\to\infty} z^2 c_2(z)=\ds{\frac{1-p^2}{4}}$,\\
 $\Delta_{\infty}=\sqrt{(1-a_{\infty})^2-4b_{\infty}}=\pm p$, and $a_i=\lim_{z\to z_i}(z-z_i)c_1(z)=\ds{\frac{1}{2}}$,\\
  $b_i=\lim_{z\to z_i}(z-z_i)^2c_2(z)=0$, $\Delta_i=\sqrt{(1-1/2)^2-4.(0)}=\ds{\frac{1}{2}}$.\\
  Then, we get $t_i=-2\cos (\ds{\frac{\pi}{2}})=0$, and $t_{\infty}=2\cos (\pi\Delta_{\infty})=2\cos (\pi p)$. The values of $\Delta_i$  correspond to points $z_i$, for $i=1,\, \dots , 4$.  

The first result of our exposition is the following proposition:
\begin{prop}
\label{2Dmain}
Let   $p=\sqrt{1+\frac{4F}{E}}$ are rational numbers, then the
system (\ref{1.2}) has no an additional meromorphic first integral when   $N(p)\ge 4$.

\end{prop}

{\bf Proof:}
 We obtain that $t_{\infty}\notin  \mathbb{Q}[t_1,\, t_2,\, t_3,\, t_4]=\mathbb{Q}$, i. e. $t_{\infty}$ is transcendental over $\mathbb{Q}$. (See \cite{Chur} for details.)

  This finishes  the proof of $b)$ by Theorem \ref{main2D_1}.

\subsection{Hamiltonian systems with homogeneous potentials.}
The study of integrability of two-dimensional systems with homogeneous potentials of different degrees turns out to be a very interesting and important problem for many researchers (see \cite{DarbouxPoint1}, \cite{DarbouxPoint}, \cite{Jarno1}, \cite{Jarno}, \cite{NakYo}, \cite{Yo}, \cite{LV}, \cite{CDMP}, \cite{CK}, \cite{LMV}, \cite{WMP}, \cite{OV}). In this subsection, we will investigate how integrable potentials interact with each other, i.e. when the sum of several integrable potentials is again an integrable potential.

To proceed with the study of integrability in the two-dimensional problem, another different approach is needed, which is well developed in \cite {DarbouxPoint1}  and \cite {DarbouxPoint}. The problem reduces to Hamiltonian system with homogeneous potentials and has a well-developed scheme for research. Let explain what this method is and how we can use it for a limited number of integrable cases. We consider a Hamiltonian system
 \begin{equation}
 \label{HomPot}
H=\frac{1}{2}\left(p_r^2+p_z^2\right)+V(r,\,z) ,
  \end{equation}
with potential 
  \begin{equation*}
  V(r,\,z)=V_{min}(r,\,z)+\dots+V_{med}(r,\,z)+\dots + V_{max}(r,\,z),
  \end{equation*}
  which is a sum of homogeneous potentials. Here with $V_{min}(r,\,z)$ (respectively $V_{med}(r,\,z)$, and $V_{max}(r,\,z)$ ) we mean the smallest (some intermediate and largest ) possible degree of homogeneous part of $V(r,\,z)$. Potential $V(r,\,z)$ is called integrable if its corresponding Hamiltonian system (\ref{HomPot}) is integrable. As it is noted in \cite{Jarno}, \cite{NakYo} and \cite {DarbouxPoint}, if $V(r,\,z)$ is an integrable potential, then $V_{min}(r,\,z)$, $V_{med}(r,\,z)$ and $V_{max}(r,\,z)$ are also integrable (for each $V_{med}(r,\,z)$). This observation gives us a chance to reduce the number of free parameters. In fact, we study for integrability the sum of proportional to integrable homogeneous potentials of degrees two, three and four. We start from the sum of integrable potentials of the second and third degrees, and then we will add the integrable potentials of the fourth degree to them. In our cases, at least for the second (it is integrable) and third degrees, the potentials are not full as in the general case, which simplifies the research a little, for the fourth degree, however, we are not so lucky. 
  In our case we have $V_{min}(r,\,z)=Ar^2+Bz^2$,  $V_{med}(r,\,z)=Cz^3+Dr^2z$, $V_{max}(r,\,z)=Ez^4+Fr^2z^2+Gr^4$ and since $V_{min}(r,\,z)$ is integrable, then the possible integrable cases are those for which $V_{max}(r,\,z)$ is integrable. These potentials are fully investigated in \cite{DarbouxPoint} and using the notation in this paper the possible integrable cases in $V_{max}(r,\,z)$ are: $V_0$, $V_1$, $V_3$, $V_4$, $V_5$ and $V_6$. Let us consider these cases in the context of the fourth degree potential in our task.

   Now we have the case {\it c.0} $V_0:=0$. The integrable homogeneous potentials of degree 3 in our case are $z^3+3r^2z$, $2z^3+r^2z$, and $16z^3+3r^2z$ (see \cite{DarbouxPoint1} for details). The Hamiltonians in this cases  are
 \begin{eqnarray*}
 H & =& \frac{1}{2}(p_r^2+p_z^2)+Ar^2+Bz^2+z^3+3r^2z,\\
  H & =& \frac{1}{2}(p_r^2+p_z^2)+Ar^2+Bz^2+2z^3+r^2z,\\
  H & =& \frac{1}{2}(p_r^2+p_z^2)+Ar^2+Bz^2+16z^3+3r^2z,
 \end{eqnarray*}
the  equations of motion  for first Hamiltonian are:
\begin{align}
\label{A0}
\dot r  &=  p_r,\, \dot p_r  =  -(2Ar+6rz),\nonumber  \\
\dot z  &=  p_z,\, \dot p_z  =  -(2Bz+3r^2+3z^2).
\end{align} 
We find a partial solution for (\ref{A0}), we put
$r = p_r = 0$ in (\ref{A0}) and we obtain
$$
\ddot z=-(2Bz+3z^2),
$$
and 
\begin{equation}
\label{A0.1}
{\dot z}^2=-2(z^3+Bz^2+h),
\end{equation}
where $h$ is a  constant. Next we change  the variable $t$ to $z$ in NVE and we obtain
 \begin{equation}
\label{NVE.A0}
 \xi_{11}''  +\frac{3z^2+2Bz}{2(z^3+Bz^2+h)}\xi_{11}'  -\frac{A+3z}{z^3+Bz^2+h}\xi_{11}=0. 
\end{equation}
The  equation (\ref{NVE.A0}) is Fuchsian (Lam\'e equation) with regular singular points at the roots $z_i$ ($z_i\ne z_j$, for $i\ne j$) of $z^3+Bz^2+h=0$ and at infinity. 
 The  indicial equations to the singular points $z_i$ and $z=\infty$ are
 $\lambda^2-\frac{1}{2}\lambda=0$, roots area
$\lambda_1=0$, $\lambda_2=\frac{1}{2}$,  for $z_i$, second one is
$2\rho^2-\rho-6 =0$, with roots
$\rho_1=2$, and $\rho_2=-\frac{3}{2}$ for
$z=\infty$.
At infinity (changing $x:=\frac{1}{z}$) equation (\ref{NVE.A0}) assumes the following form:
\begin{equation}
\label{NVE.A0Infty}
 \frac{d^2\xi_{11}}{dx^2}  +\frac{4hx^3+2Bx+1}{2x(hx^3+Bx+1)}\frac{d\xi_{11}}{dx}  -\frac{Ax+3}{x^2(hx^3+Bx+1)}\xi_{11}=0,
\end{equation}
 with  a local solutions at $x=0$
\begin{eqnarray*}
\xi_{11}^{(1)}(x) & = & x^{-3/2}\left(1+(\frac{-2A}{5}+\frac{9B}{10})x+(\frac{2}{15}A^2-\frac{1}{3}AB+\frac{3}{40}B)x^2\dots\right),\\
\xi_{11}^{(2)}(x) & = &x^2\left(1+(\frac{2A}{9}-\frac{8B}{9})x+(\frac{2}{99}A^2-\frac{26}{99}AB+\frac{8}{11}B)x^2\dots\right).
  \end{eqnarray*}
  Now let us repeat the procedure for the second equation in $VE_1$, we get
  \begin{equation}
\label{NVE.A0_2}
 \xi_{12}''  +\frac{3z^2+2Bz}{2(z^3+Bz^2+h)}\xi_{12}'  -\frac{B+3z}{z^3+Bz^2+h}\xi_{12}=0, 
\end{equation}
  \begin{equation}
\label{NVE.A0_2Infty}
 \frac{d^2\xi_{12}}{dx^2}  +\frac{4hx^3+2Bx+1}{2x(hx^3+Bx+1)}\frac{d\xi_{12}}{dx}  -\frac{Bx+3}{x^2(hx^3+Bx+1)}\xi_{12}=0,
\end{equation}
with  a local solutions at $x=0$
\begin{eqnarray*}
\xi_{12}^{(1)}(x) & = & x^{-3/2}\left(1+\frac{B}{2}x-\frac{1}{8}B^2x^2\dots\right),\\
\xi_{12}^{(2)}(x) & = &x^2\left(1-\frac{2B}{3}x+\frac{16}{33}B^2x^2\dots\right).
  \end{eqnarray*}
 For $VE_2$ we do the same transformation of variables $t\rightarrow z\rightarrow x$ and we obtain 
  
  \begin{eqnarray*}
\label{VE2.A0_2Infty}
 \frac{d^2\xi_{11}}{dx^2}  +\frac{4hx^3+2Bx+1}{2x(hx^3+Bx+1)}\frac{d\xi_{11}}{dx}  -\frac{Ax+3}{x^2(hx^3+Bx+1)}\xi_{11}=\frac{\xi_{11}\xi_{12}}{2x(hx^3+Bx+1)},\\
 \frac{d^2\xi_{12}}{dx^2}  +\frac{4hx^3+2Bx+1}{2x(hx^3+Bx+1)}\frac{d\xi_{12}}{dx}  -\frac{Bx+3}{x^2(hx^3+Bx+1)}\xi_{12}=\frac{3\xi_{11}^2+3\xi_{12}^2}{2x(hx^3+Bx+1)}.
\end{eqnarray*} 
The necessary condition for integrability is the absence of a residue in the expressions 
$$\xi_{11}.\frac{\xi_{11}\xi_{12}}{\sqrt{2x(hx^3+Bx+1)}}$$
 and 
 $$\xi_{12}.\frac{3\xi_{11}^2+3\xi_{12}^2}{\sqrt{2x(hx^3+Bx+1)}},$$ 
 which is 
 $$\frac{32(A-B).B.(A^2-\frac{17}{8} .AB-\frac{9}{512} . B^2)}{225},$$
 for solutions $\xi_{11}^{(1)}(x)$ and $\xi_{12}^{(1)}(x)$. 
 Then we get three cases of possible integrability: $A=B$, $B=0$ and $A^2-\frac{17}{8} .AB-\frac{9}{512} . B^2=0$. By using third variations it turns out that the second and third cases are non-integrable at $A\ne 0$.

Let us look at the main fragments of the proof in the case when $B=0$: The Hamiltonian and equations of motion are 
$$H=\frac{1}{2}(p_r^2+p_z^2)+Ar^2+z^3+3r^2z,$$
\begin{align*}
\dot r  &=  p_r,\, \dot p_r  =  -(2Ar+6rz),\nonumber  \\
\dot z  &=  p_z,\, \dot p_z  =  -(3z^2+3r^2).
\end{align*}
For partial solution in this case we  can choose $p_r=r=0$ and we obtain $\dot{z}^2=-2z^3-2h$ ($h=const\ne 0$). After changing the variables $z=-2\tilde{z}$ we get $\dot{\tilde{z}}^2=4\tilde{z}^3-\frac{h}{2}$ which can be written in the form $\tilde{z}=\wp(t,0,\frac{h}{2})$ using the Weierstrass p-function. For the variational equations $VE1$, $VE2$ and $VE3$  are obtained as follows
\begin{align*}
 \ddot {\xi_{11}} =&  \left(12\wp(t,0,\frac{h}{2})-2A\right)\xi_{11},\nonumber  \\
 \ddot {\xi_{12}} = & \left(12\wp(t,0,\frac{h}{2})\right)\xi_{12},
\end{align*} 
\begin{align*}
 \ddot {\xi_{21}} =&  \left(12\wp(t,0,\frac{h}{2})-2A\right)\xi_{21}-6\xi_{11}\xi_{12},\nonumber  \\
 \ddot {\xi_{22}} = & \left(12\wp(t,0,\frac{h}{2})\right)\xi_{22}-3\left((\xi_{11})^2+(\xi_{12})^2\right),
\end{align*} 
and
 \begin{align*}
 \ddot {\xi_{31}} = & \left(12\wp(t,0,\frac{h}{2})-2A\right)\xi_{31}-6\left(\xi_{11}\xi_{22}+\xi_{12}\xi_{21}\right),\nonumber  \\
 \ddot {\xi_{32}} =&  \left(12\wp(t,0,\frac{h}{2})\right)\xi_{32}-6\left(\xi_{11}\xi_{21}+\xi_{12}\xi_{22}\right).
\end{align*} 
For solutions of $VE1$ around $t=0$ we have
\begin{eqnarray*}
\xi_{11}^{(1)}(t)  = & t^4-\frac{A}{9}t^5\dots,\\
\xi_{11}^{(2)}(t)  = &t^{-3}+\frac{A}{5}t^{-1}\dots,
  \end{eqnarray*}
\begin{eqnarray*}
\xi_{12}^{(1)}(t)  = & t^4-\frac{h}{364}t^{10}\dots,\\
\xi_{12}^{(2)}(t)  = &t^{-3}+\frac{h}{28}t^{3}\dots.
  \end{eqnarray*}
Let us use notation
\begin{align*}
K_2^{(1)} = & -6\left(\xi_{11}\xi_{12}\right), \\
K_2^{(2)} = & -3\left((\xi_{11})^2+(\xi_{12})^2\right), \\
K_3^{(1)} = & -6\left(\xi_{11}\xi_{22}+\xi_{12}\xi_{21}\right), \\
K_3^{(2)} = & -6\left(\xi_{11}\xi_{21}+\xi_{12}\xi_{22}\right),
\end{align*} 
and
 \begin{align*}
f_2 = & \left(0,\, K_2^{(1)},\, 0,\, K_2^{(2)} \right)^T, \\
f_3 = & \left(0,\, K_3^{(1)},\, 0,\, K_3^{(2)} \right)^T.
\end{align*} 
Without losing a community we can assume that
$\xi_{11}^{(1)}\dot{\xi_{11}^{(2)}}-\xi_{11}^{(2)}\dot{\xi_{11}^{(1)}}=1$ and $\xi_{12}^{(1)}\dot{\xi_{12}^{(2)}}-\xi_{12}^{(2)}\dot{\xi_{12}^{(1)}}=1$.
Then the fundamental matrix of $VE1$ and its inverse are
\begin{equation}
\label{X}
X (t) =
\begin{pmatrix}
 \xi_{11} ^{(1)} & \xi_{11} ^{(2)} & 0 & 0\\
 \dot{\xi}_{11} ^{(1)} & \dot{\xi}_{11} ^{(2)}   & 0 & 0 \\
0 & 0 & \xi_{12} ^{(1)} & \xi_{12} ^{(2)}\\
0 & 0 & \dot{\xi}_{12} ^{(1)} & \dot{\xi}_{12} ^{(2)}
  \end{pmatrix},
  \end{equation}

 \begin{equation}
\label{X^1}
X ^{-1}(t) =
\begin{pmatrix}
 \dot{\xi}_{11} ^{(2)} & -\xi_{11} ^{(2)} & 0 & 0\\
 -\dot{\xi}_{11} ^{(1)} & \xi_{11} ^{(1)}   & 0 & 0 \\
0 & 0 & \dot{\xi}_{12} ^{(2)} & -\xi_{12} ^{(2)}\\
0 & 0 & -\dot{\xi}_{12} ^{(1)} & \xi_{12} ^{(2)}
  \end{pmatrix}.
  \end{equation}
We look at expression $X^{-1}f_2$ for a non-zero residue - there is none in all possibilities, we do the same for expression $X^{-1}f_3$, there is a non-zero residue for $A\ne 0$ for the expression $\xi_{11}^{(1)}\left(\xi_{11}^{(2)}\xi_{22}^{(1)}+\xi_{12}^{(2)}\xi_{21}^{(1)}\right)$ of $X^{-1}f_3$.

For the second Hamiltonian in this case, we have not integrability constraints are obtained under the second variations.

For the third Hamiltonian in this case, we consistently obtain
\begin{align*}
\dot r  &=  p_r,\, \dot p_r  =  -(2Ar+6rz),\nonumber  \\
\dot z  &=  p_z,\, \dot p_z  =  -(2Bz+48z^2+3r^2)
\end{align*} 
The partial solution in this case is $r = p_r = 0$ and $\ddot z=-(2Bz+48z^2)$, ${\dot z}^2=-2(16z^3+Bz^2+h)$, $NVE$ in $\infty$ is 
\begin{equation}
\label{NVE.A1Infty_1}
 \frac{d^2\xi_{11}}{dx^2}  +\frac{2hx^3+8}{x(hx^3+Bx+16)}\frac{d\xi_{11}}{dx}  -\frac{Ax+48}{x^2(hx^3+Bx+16)}\xi_{11}=0.
\end{equation}
The   a local solutions  of last equation at $x=0$ are
\begin{eqnarray*}
\xi_{11}^{(1)}(x) & = & x^{3/4}\left(1+(\frac{A}{32}+\frac{3B}{512})x+\dots\right),\\
\xi_{11}^{(2)}(x) & = &x^{-1/4}\left(\ln(x).(16A-5B)(\frac{1}{16}x+\dots)+\dots\right).
  \end{eqnarray*}
From Fuchs's Theorem, if we have a logarithm in one of solutions of linear differential equation of the second order then
\begin{equation*}
\xi_{11}^{(2)}(x)=\xi_{11}^{(1)}(x)\ln (x) +h(x), 
 \end{equation*}
 where $h(x)$ is a holomorphic function (in this case, or meromorphic in in the general case).
  We denote with $\Gamma$ the Riemann surface $\left\{(x,y):\,  y^2=2hx^4+2Bx^2+16x\right\}$. The field of coefficients of (\ref{NVE.A1Infty_1}) is K:=$\mathcal{M}(\Gamma)$ (meromorphic functions on $\Gamma$). Let us also denote by $L:=K[\xi_{11}^{(1)}, \xi_{11}^{(2)}]$ the Picard-Vesiot extension of $K$. Let us formulate the following 
  \begin{lm}
  \label{lemma1}
  In the above notations for (\ref{NVE.A1Infty_1}), if $\sigma\in Gal(L/K)$, then $\sigma(\ln (x))=\delta\ln(x)+\gamma$, for $\delta, \,\gamma\in\mathbb{C} $.
  \end{lm}
  
  {\bf Proof:}
  We have
   $$\ds{\frac{d\sigma(\ln( x))}{dx}=\sigma(\frac{d\ln (x)}{dx})=\sigma(\frac{1}{x})},$$
  $\sigma ' .\ds{\frac{1}{x}}=\sigma(\ds{\frac{1}{x}})$, let $x=e^t$, 
  $\ds{\frac{d}{dt}}\sigma . e^{-t}=\sigma(e^{-t})=\delta e^{-t}$, $\delta\in \mathbb{C}$, then we have $\ds{\frac{d\sigma}{dt}}=\delta$, and $\sigma(t)=\delta t+\gamma$, 
  then $\sigma(\ln (x))=\delta\ln(x)+\gamma$, for $\delta, \,\gamma\in\mathbb{C} $.
  
  For the local Galois group of (\ref{NVE.A1Infty_1}) we have:
 $\sigma\in Gal(L/K)$, $\sigma(\xi_{11}^{(1)})=\ds{\frac{1}{\delta}}\xi_{11}^{(1)}$, $\delta \in \mathbb{C^*}$, 
 \begin{eqnarray*}
 \sigma(\xi_{11}^{(2)}) & = & \frac{1}{\delta}\xi_{11}^{(1)}\sigma(\ln(x))+h(x)\\
 &= & \frac{1}{\delta}\xi_{11}^{(1)}(\delta\ln(x)+\gamma)+h(x)\\
  &= & \frac{1}{\delta}\xi_{11}^{(1)}(\delta\ln(x))+h(x)+\frac{\gamma}{\delta}\xi_{11}^{(1)}\\
  &= & \xi_{11}^{(2)}+\frac{\gamma}{\delta}\xi_{11}^{(1)} .
\end{eqnarray*}
Let $X(x)$ is the fundamental matrix of (\ref{NVE.A1Infty_1}), then
$\sigma (X(x))=X(x)R $, and $R= \left\{\begin{pmatrix}
\frac{1}{\delta} & 0  \\
\frac{\gamma}{\delta} & 1
\end{pmatrix},  \,\delta\ne0\right\}$ - solvable but non-commutative. It is not difficult to observe that the condition for non-commutativity of the local Galois group is equivalent to the existence of a logarithm in one of the solutions of (\ref{NVE.A1Infty_1}). In the present case, the condition for the existence of a logarithm is $16A\ne 5B$. In order to prove the non-commutativity of the global Galois group of (\ref{NVE.A1Infty_1}), we need to study the local group near of the regular singular point $z_1$ ($z_1$ is the root of $z^4+Cz^3+Bz^2+h=0$). We know the roots of the indicative equation($\lambda_1=0$ and $\lambda_2=\frac{1}{2}$) and this allows us to conclude that one local solution is in the field of constants and the other is not. From this we obtain that the local group of Galois consists of matrices of type $\left\{\begin{pmatrix}1 & 0  \\
\mu & 1
\end{pmatrix},  \,\mu \ne 0\right\}$, which in the general case do not commute with those already found in infinity  $ \left\{\begin{pmatrix}
\frac{1}{\delta} & 0  \\
\frac{\gamma}{\delta} & 1
\end{pmatrix},  \,\delta\ne0\right\}$ at $16A\ne 5B$. (See \cite{MR1} for details.)

 The case {\it c.0)} is proved.

Unfortunately, the proof of points {\it c)} to the theorem \ref{main2D_1} is quite extensive (a full consideration of each case deserves a separate study), so we will mention that all methods of studying (second and third variations and existence of logarithm) were already used to prove the case {\it c.0)}.

\section{Remarks and Comments}
Given the above statement, we can easily conclude that the integrable cases in \cite{Comment I} are exceptions to the already obtained non-integrability conditions in Theorem \ref{main2D_1} (the Cases 1, 2 and 3 in \cite{Comment I} correspond to c.3), c.8) and c.4) respectively).
Now, to clarify, where do such significant differences in the results of publications \cite{Comment I} and \cite{G} come from? After a careful analysis, it turns out that there are three reasons: the first is that in \cite{G} a very special case h = 0 (zero energy level) is considered, which, although quite interesting as a study, distorts the non-integrability result ;

The second reason for the differences is the ignoring in \cite{G} of the logarithmic terms in the local solutions - it turns out that they play an important role in the study of the Galois Group of the variational differential equations under consideration;

The third reason is the behavior of sums of different integrable homogeneous potentials - It turns out that not always a sum of homogeneous integrable potentials is an integrable potential.
I'd also like to note that studying cases with recurrent roots has proven (thanks to a large number of computations) to be quite a large challenge, and we could frame them as open problems. 

From what is shown in \cite{G} and the examples in \cite{Comment I}, we can also conclude something that I think is important, namely that zero-level ($h=0$) non-integrability does not necessarily imply non-integrability in the general case ($h\ne 0$).

It remains to notice that the first integrals found in \cite{Trapped} and \cite{Comment I} are too limited in number against the background of the integrals that could be found (if we seek for  them in polynomial form, because there exists a technique for this, although not so simple). The searching for such integrals is the next open question of this study.

 \section*{ Acknowledgments }

The author would like to express special thanks to Dr. Idriss Elfakkousy for useful discussions.
The author has been partially supported by  grant  80-10-53/10.05.2022 from Sofia University ``St. Kliment Ohridski''.

\end{document}